

%




\input amstex
\documentstyle{amsppt}
\magnification=\magstep1
\NoBlackBoxes
\settabs 10 \columns

\def\rarrow{\rightarrow}
\def\e{\epsilon }
\def\ms{\medskip}
\def\bs{\bigskip}

\def\ss#1{_{\lower3pt\hbox{$\scriptstyle #1$}}}
\def\nm#1{\left\|#1\right\|}

\def\o{\over}
\def\Ball{\hbox{Ball }}
\def\twosum{$2$-summing}
\def\twosump{$2$-summing property}
\def\Rn#1{\hbox{{\it I\kern -0.25emR}$\sp{\,{#1}}$}}
\def\sm#1{\sum_1^#1}
\def\seq#1#2{#1_1,\dots,#1_#2}
\def\qed {\vrule height6pt  width6pt depth0pt}
\def\scpr#1{\langle#1\rangle}
\def\cal{\Cal}
\def\span{\hbox{span }}
\def\R{\hbox{{\it I\kern -0.25emR}}}

\def\T{{\Bbb T}}

\def\C{{\Bbb C}}
\def\e{\epsilon }
\def\a{\alpha }
\def\b{\beta }
\def\d{\delta }
\def\g{\gamma }

\def\lam{\lambda}
\def\longlongrightarrow{\relbar\joinrel\relbar\joinrel
\relbar\joinrel\rightarrow}
\topmatter

\title Banach spaces with the $2$-summing property
\endtitle

\author A. Arias, T. Figiel, W.~B. Johnson and G. Schechtman
\endauthor

\thanks
The first author was supported as Young Investigator, NSF DMS 89-21369; and
NSF DMS 93-21369;
the second and fourth authors were Workshop participants, NSF \#DMS 89-21369;
the third author was supported by NSF \#DMS 90-03550 and 93-06376;
the third and fourth authors were supported by the U.S.-Israel
BSF
\endthanks

\address A. Arias
\newline Division of Mathematics, Computer Science and Statistics, 
The University of Texas at San Antonio, San
Antonio, TX 78249, U.S.A.\endaddress
\email arias\@ringer.cs.utsa.edu\endemail

\address T.~Figiel
\newline Institute of Mathematics, Polish Academy of Sciences, Gdansk,
Poland \endaddress
\email tfigiel\@impan.impan.gov.pl\endemail

\address W.~B.~Johnson
\newline Department of Mathematics, Texas A\&M University, College Station, TX
77843, U.S.A. \endaddress
\email johnson\@math.tamu.edu\endemail

\address G.~Schechtman
\newline Department of Theoretical Mathematics, The Weizmann
Institute of Science, Rehovot, Israel\endaddress
\email mtschech\@weizmann.weizmann.ac.il\endemail

\abstract A Banach space $X$ has the $2$-summing property if
the norm of every linear operator from $X$ to a Hilbert space is
equal to the $2$-summing norm of the operator.
Up to a point, the theory of spaces which have this property
is independent of the scalar field: the property is
self-dual and any space with the property  is a finite dimensional space of
maximal distance to the Hilbert space of the same dimension.
In the case of real scalars  only  the real line and real
$\ell_\infty^2$ have the $2$-summing property. In the complex case there are
more
examples; e.g., all subspaces of complex $\ell_\infty^3$ and their duals.
\endabstract


\endtopmatter

\document
\heading 0. Introduction:\endheading

Some important classical Banach spaces; in
particular, $C(K)$ spaces, $L_1$ spaces, the disk algebra; as well as some
other spaces (such as quotients of $L_1$ spaces by reflexive subspaces
[K], [Pi]), have the property that every (bounded, linear) operator from
the space into a Hilbert space is $2$-summing. (Later we review equivalent
formulations of the definition of $2$-summing operator. Here we mention
only that an operator  $T : X \rarrow \ell_2$ is $2$-summing provided that
for all operators  $u : \ell_2 \rarrow X$ the composition $Tu$ is a
Hilbert-Schmidt operator; moreover, the $2$-summing norm $\pi_2(T)$ of
$T$ is the supremum of the Hilbert-Schmidt norm of $Tu$ as $u$ ranges
over all norm one operators $u : \ell_2 \rarrow X$.)  In this paper we
investigate the isometric version of this property: say  that a Banach
space $X$ has the {\it $2$-summing property\/} provided that
$\pi_2(T) = \nm{T}$ for all operators  $T : X \rarrow \ell_2$.

While the $2$-summing property is a purely Banach space concept and our
investigation lies purely in the realm of Banach space theory, part of
the motivation for studying the $2$-summing property comes from operator
spaces. In [Pa], Paulsen defines for a Banach space $X$ the parameter
$\a(X)$ to be the supremum of the completely bounded norm of $T$ as $T$
ranges over all norm one operators from $X$ into the space $B(\ell_2)$
of all bounded linear operators on $\ell_2$ and asks which spaces $X$
have the property that $\a(X)=1$. Paulsen's
problem and study of $\a(X)$ is motivated by old results of von Neumann,
Sz.-Nagy,  Arveson, and Parrott as well as more recent research of Misra
and Sastry.   The connection between Paulsen's problem and the present
paper is Blecher's result [B] that $\a(X)=1$  implies that $X$ has
the  $2$-summing property. Another connection is through the {\it
property (P)\/} introduced by Bagchi and Misra [BM], which Pisier
noticed is equivalent to the $2$-summing property.   However, since we
shall not investigate here $\a(X)$ or property (P) directly and do not
require results from  operator theory, we refer the interested reader to
[Pa] and [BM] for definitions, history, and references.  On the other
hand,  since the topic we treat here is relevant for operator theorists,
we repeat standard background in Banach space theory used herein for
their benefit.

In [Pa] Paulsen asks whether $\a(X)=1$ only if $X$ is a one or two
dimensional $C(K)$ or $L_1$ space; in other words, ignoring the trivial
one dimensional case, whether $\a(X)=1$ implies that $X$ is two
dimensional, and among two dimensional spaces, whether only
$\ell_{\infty}^2$ and $\ell_1^2$ satisfy this identity. He proves that
$\a(X)=1$ implies that $\dim(X)$ is at most $4$, that $\a(X)=\a(X^*)$,
and he gives another proof of Haagerup's theorem that
$\a(\ell_{\infty}^2)=1$.  Paulsen, interested in operator theory, is
referring to {\sl complex\/} Banach spaces, so $\ell_{\infty}^2$ is not
the same space as $\ell_1^2$.

From the  point of view of Banach space theory, it is natural to ask
which Banach spaces have the $2$-summing property both in the real and
the complex cases, and here we investigate both questions. Up to a
point, the theory is independent of the scalar field: In section 2 we
show that the  \twosump\ is self-dual, that only spaces of sufficiently
small (finite) dimension can have the property, and that a space  with
the property is a maximal distance space--that is, it has maximal
Banach-Mazur distance to the Hilbert space of the same dimension. The
main result in section 2, Proposition 2.6, gives a useful condition for
checking whether a space has the \twosump\ which takes a particularly
simple form when the space is $2$-dimensional (Corollary 2.7.a).

The analysis in section 3 yields that the situation is very simple in
the case of real scalars; namely,  $\R$ and $\ell_{\infty}^2$ are the
only spaces which have the \twosump. Two ingredients for proving this are
Proposition  3.1, which says that there are many norm one operators from
real $\ell_1^3$ into $\ell_2^2$ which have \twosum\ norm larger than
one, and a geometrical argument, which together with a recent lemma of
Maurey implies that a maximal distance real space of dimension at least
three has a two dimensional quotient whose unit ball is a regular hexagon.

The complex case is {\sl a priori\/} more complicated, since
$\ell_{\infty}^2$ and $\ell_1^2$ both have the \twosump\ but are not
isometrically isomorphic.  In fact, in section 4 we show that there are
many other examples of complex spaces which have the \twosump; in
particular, $\ell_{\infty}^3$ and all of its subspaces.  The simplest way
to prove that these spaces have the \twosump\ is to apply Proposition 2.6,
but we also give direct proofs for $\ell_{\infty}^3$ in section 4
and for its two dimensional subspaces in the appendix. The case of
$\ell_{\infty}^3$ itself reduces via a simple but slightly strange
``abstract nonsense" argument to a calculus lemma, which, while easy,
does not look familiar.  (In [BM] the authors give an argument that
$\ell_{\infty}^3$ satisfies their property (P) which uses a variation of
the calculus lemma but replaces the ``abstract nonsense" with a
reduction to self-adjoint matrices.)   We also give in Proposition 4.5
an inequality which is equivalent to the assertion that all two
dimensional subspaces of complex $\ell_{\infty}^3$ have the \twosump.
While we do not see a simple direct proof of this inequality, we give a
very simple proof of a weaker inequality which is equivalent to the
assertion that every two dimensional subspace  of complex
$\ell_{\infty}^3$ is of maximal distance.

In section 5 we make some additional observations.
\bs

\head 1. Preliminaries.\endhead
Standard Banach space theory language and results
can be found in  [LT1], [LT2],  while basic results in the  local theory
of Banach spaces are contained in [T-J2].  However, we recall here that
part of the theory and language which we think may not be well known to
specialists in operator theory.

Spaces are always Banach spaces, and subspaces are assumed to be closed.
Operators are always bounded and linear.  The [Banach-Mazur] distance
between spaces $X$ and $Y$ is the isomorphism constant, defined as the
infimum of $\nm{T}\nm{T^{-1}}$ as $T$ runs over all invertible operators
from $X$ onto $Y$. The closed unit ball of $X$ is denoted by $\Ball(X)$.
``Local theory" is loosely defined as the study of
properties of infinite dimensional spaces which depend only on their
finite dimensional spaces, as well as the study of numerical parameters
associated with finite dimensional spaces. Basic for our study and most
other investigations in local theory is the fact (see [T-J2,~p.~54]) that
the distance from an $n$-dimensional space to $\ell_2^n$ is at most
$\sqrt{n}$.  One  proves this by using the following consequence of F.
John's theorem ([T-J2,~p.~123]): If $X$ is  $n$-dimensional and $\cal E$
is the ellipsoid of minimal volume containing $\Ball(X)$, then $n^{-{1\o
2}}{\cal E} \subset \Ball(X)$. This statement perhaps should be
elaborated:  Since $\dim(X) < \infty$, we can regard $X$ as $\R^n$ or
$\C^n$ with some norm. Among all norm-increasing operators $u$ from
$\ell_2^n$ into $X$, there is by compactness one which minimizes the
volume of $u(\Ball(\ell_2^n))$; the distance assertion says that
$\nm{u^{-1}}\le \sqrt{n}$.    Alternatively, if one chooses from among
all norm one operators from $X$ into $\ell_2^n$ one which maximizes the
volume of the image of $\Ball(X)$, then the norm of the inverse of this
operator is at most $\sqrt{n}$.   If complex $\ell_2^n$ is considered as
a real space, then it is isometrically isomorphic to real $\ell_2^{2n}$.
Thus the distance statement for {\sl complex\/} spaces says that a
complex space of dimension $n$, when considered as a {\sl real\/} space
of dimension $2n$,  has (real) distance to (real) $\ell_2^{2n}$ at most
$\sqrt{n}$.

Actually, we need more than just the distance consequence of John's
theorem. The theorem itself [T-J2,~p.~122] says that if $\cal E$ is the
ellipsoid of minimal volume containing $\Ball(X)$, then there exist
points of contact  $\seq{y}{m}$ between the unit sphere of $X$ and the
boundary of $\cal E$, and there exist positive real numbers
$\seq{\mu}{m}$ summing to $\dim X$ so that for each $x$ in $X$,
$x=\sum_{i=1}^m \mu_i\scpr{x,y_i}y_i$, where ``$\scpr{\cdot,\cdot}$" is
the scalar product which generates the ellipsoid $\cal E$.  The existence
of many contact points between $\Ball(X)$ and $\cal E$ is important for
the proof of Theorem 3.3.

The dual concept to minimal volume ellipsoid is maximal volume
ellipsoid.  More precisely, an $n$-dimensional space can be regarded as
$\R^n$ or $\C^n$ under some norm $\nm{\cdot}$ in such a way that
$\Ball(X)\subset \cal E$, where $\cal E$ is the usual Euclidean ball and
is also the ellipsoid of minimal volume containing $\Ball(X)$.  Then
$X^*$ is naturally represented as $\R^n$ or $\C^n$ under some norm, and
the action of $X^*$ on $X$ is given by the usual inner product.  Then
$\cal E$ is the ellipsoid of maximal volume contained in $\Ball(X^*)$.

John's theorem gives many points of contact between $\Ball(X)$ and the
boundary of the ellipsoid of minimal volume containing $\Ball(X)$, and
many points of contact between the boundary of $\Ball(X)$ and the
ellipsoid of maximal volume contained in $\Ball(X)$. It is a nuisance
that these two ellipsoids are not generally homothetic; however, the
situation is better when $X$ has the \twosump:

\proclaim {Lemma 1.1}  Assume that the real or complex $n$-dimensional
space $X$ has the \twosump \ and let ${\Cal E}_1$ be the ellipsoid of
minimal volume containing $\Ball(X)$.  Then $n^{-{1\o 2}}{\Cal E}_1$ is
the ellipsoid of maximal volume contained in $\Ball(X)$.
\endproclaim

{\bf Proof:} Let ${\Cal E}_2$ be the ellipsoid of maximal volume
contained in  $\Ball(X)$ and for $i=1,2$ let $|\cdot|_i$ be the
Euclidean norm which has for its unit ball ${\Cal E}_i$.   Let $u_1$ be
the formal identity map from $X$ to the Euclidean space $(X,|\cdot|_1)$,
$u_2$ the identity map from  $(X,|\cdot|_2)$  to  $X$, and let
$\seq{\lam}{n}$ be the $s$-numbers of the Hilbert space operator
$u_1u_2$ (i.e., the square roots of the eigenvalues of
$(u_1u_2)^*u_1u_2$).
 Since $\pi_2(u_1)=\|u_1\|=1$ and $\|u_2\|=1$ we have that
$\pi_2(u_1u_2)\leq1$. This implies that
$$
|\lambda_1|^2+\cdots+|\lambda_n|^2\leq1.\leqno{(1.1)}
$$
On the other hand  $\hbox{vol}({\cal E}_2)\ge\hbox{vol}({\cal E}_1/\sqrt{n})$
because $n^{-{1\o 2}}{\cal E}_1 \subset \Ball(X)$, so that in the case
of real scalars we get,
$$
\lambda_1\lambda_2\cdots\lambda_n\geq\biggl({1\over\sqrt{n}}\biggr)^n,
\leqno{(1.2.R)}
$$
and in the complex case
$$
\lambda_1^2\lambda_2^2\cdots\lambda_n^2\geq\biggl({1\over\sqrt{n}}\biggr)^{2n}.
\leqno{(1.2.C)}
$$
The only way that (1.1) and (1.2) are true is if
$\lambda_1=\lambda_2=\cdots =\lambda_n=1/\sqrt{n}$. But this implies that
${\cal E}_2={\cal E}_1/\sqrt{n}$.  \hfill\qed

\ms

{\bf Remark 1.2.}  B.~Maurey has proved a far reaching generalization of
Lemma~1.1; namely, that if a space $X$ does not have a unique (up to
homothety) distance ellipsoid, then there is a subspace which has the
same distance to a Hilbert space as the whole space and   which has a
unique  distance ellipsoid.  This implies an unpublished result due to
Tomczak-Jaegermann which is stronger than Lemma 1.1; namely, that when the
distance is maximal, the minimal and maximal volume ellipsoids must be
homothetic.

\medbreak

Basic facts about \twosum\ operators, and, more generally, $p$-summing
operators, can be found in [LT1] and [T-J2].  The $2$-summing norm
$\pi_2(T)$ of an operator from a space $X$ to a space $Y$ is defined to
be the supremum of $\left(\sm{n} \nm{TUe_i}^2\right)^{1/2}$ where the sup
is over all norm one operators $U$ from $\ell_2^n$, $n=1,2,\dots$, into
$X$ and $\{e_i\}_{i=1}^n$ is the unit vector basis for $\ell_2^n$.  When
$Y$ is a Hilbert space, this reduces to the definition given in the first
paragraph of the introduction, and when $X$ is also a Hilbert space,
$\pi_2(T)$ is the Hilbert-Schmidt norm of $T$. Note that if  $U$  is an
operator from $\ell_2^n$ to a subspace $X$ of $\ell_\infty$, then
$\nm{U}^2 = \nm{\sum_{i=1}^n|Ue_i|^2}$ [the absolute value is
interpreted coordinatewise in $\ell_\infty$]. So if $T$ goes from $X$
into a space $Y$, $\pi_2(T)$ can be defined intrinsically by $$
\pi_2(T)^2 = \sup \{\sum_{i=1}^n\nm{Tx_i}^2 : \ \
\Bigl\Vert\sum_{i=1}^n|x_i|^2\Bigr\Vert\le 1; \ x_i \in X; \
n=1,2,3,...\}; $$ but when $Y$ is an $N$-dimensional Hilbert space, the
``sup" is already achieved for $n = N$.  (Not relevant for this paper
but worth noting is that when $Y$ is a general $N$-dimensional space,
the ``sup" is achieved for $n \le N^2$ [FLM], [T-J2, p. 141] and is estimated
up to the
multiplicative constant $\sqrt 2$ for $n = N$ [T-J1], [T-J2, p. 143].)

It is easy to see that $\pi_2$ is a complete norm on the space of all
\twosum \ operators from $X$ to $Y$ and that $\pi_2$ has the ideal
property; that is, for any defined composition $T_1T_2T_3$ of operators,
$\pi_2(T_1T_2T_3)\le \nm{T_1}\pi_2(T_2)\nm{T_3}$. The typical \twosum\
operator is the formal identity mapping $I_{\infty,2}$ from
$L_{\infty}(\Omega, \mu)$ to $L_2(\Omega,\mu)$ when $\mu$ is a finite
measure.  In this case one gets easily that
$\pi_2(I_{\infty,2})=\mu(\Omega)^{1\o 2}$.  That  such operators are
typical is a consequence of the Pietsch factorization theorem ([LT, p.~64],
[T-J2,~p.~47]), which says that if the space $X$ is isometrically
included in a $C(K)$ space, and $T:X\to Y$ is \twosum, then there is a
probability measure $\mu$ on $K$ and an operator $S$ from $L_2(K,\mu)$
into $Y$ so that $T$ is the restriction of $SI_{\infty,2}$ to $X$ and
$\nm{S}=\pi_2(T)$. That is, there is a probability measure $\nu$ on $K$ so
that for each $x$ in $X$,
$$ \nm{Tx}^2 \le \pi_2(T)^2 \int |x^*(x)|^2\, d\nu(x^*).\leqno{(1.3)}
$$
Of course, the converse to the Pietsch factorization theorem follows
from the ideal property for \twosum\ operators.

The qualitative version of Dvoretzky's theorem [T-J2,~p.~26] says that
every infinite dimensional space $X$ contains for every $n$ and $\e>0$ a
subspace whose distance to $\ell_2^n$ is less than $1+\e$.  In fact, for
a fixed $n$ and $\e$, the same conclusion is true if dim$(X)\ge
N(n,\e)$, and the known estimates for $N(n,\e)$ are rather good.

 \bs

\head 2. General results.\endhead
Here we mention some simple results about
spaces which have the \twosump, present some motivating examples and
then find a characterization of spaces with that property.  Let us say
that $X$ satisfies the $k$-dimensional \twosump\  if $\pi_2(T)=\nm{T}$
for every operator $T$ from $X$ into $\ell_2^k$.  Thus every space has
the $1$-dimensional \twosump, and $X$ has the \twosump \ if $X$ has the
$k$-dimensional \twosump\  for every positive integer $k$.  We
introduce this definition because our techniques suggest that a space
with the $2$-dimensional \twosump\ has the \twosump, but we cannot
prove this even in the case of real scalars.

Throughout this section the scalars can be either ${\R}$ or $\C$ unless
explicitly stated otherwise.

\proclaim {Proposition 2.1} \item {(a)} If $X$ has the $2$-dimensional
\twosump, then $X$ is finite dimensional.
 \item{(b)} If $X$  has the $k$-dimensional \twosump\ for some $k$,
then so does $X^*$.
\item{(c)} If $X$ has the $2$-summing property, then $X$ is a maximal
distance space.
\endproclaim

{\bf Proof:}  For (a), we use the fact that $\ell_1^m$ fails the
$2$-dimensional \twosump \ for some integer $m$.  In fact,  in the real
case, $m$ can be taken to be $3$ (Example 2.3), while in the complex
case, $m=4$ suffices (remark after example 2.3).  Alternatively,  one
can check that a quotient mapping from $\ell_1$ onto $\ell_2^2$ has
\twosum\ norm larger than one, which implies that $\ell_1^m$ fails the
$2$-dimensional \twosump \ if $m$ is sufficiently large.  So fix a norm
one operator $u$ from $\ell_1^m$ into $\ell_2^2$ for which $\pi_2(u)>1$.
By Dvoretzky's theorem, $\ell_2^2$ is almost a quotient of every
infinite dimensional space, so if $\dim X$ is sufficiently large, then
there is an operator $Q$ from $X$  into $\ell_2^2$ with $\Ball(\ell_2^2)
\subset Q[\Ball(X)]$ but $\nm{Q}<\pi_2(u)$. Pick $\seq{z}{m}$ in
$\Ball(X)$ with $Qz_i=ue_i$ for $i=1,2,\dots,m$ and define $T$ from
$\ell_1^m$ into $X$ by $Te_i = z_i$, $i=1,2,\dots,m$. Then $u=QT$ and
$\pi_2(u)\le \pi_2(Q)$ but $\|Q\| < \pi_2(u)$.

For (b), assume that $X$   has the $k$-dimensional \twosump \ and let
$T$ be any norm one operator from $X^*$ into $\ell_2^k$.  It is enough
to show that $\pi_2(Tu)\le1$ when $u$ is a norm one operator from
$\ell_2^k$ into $X^*$. This brings us back to the familiar setting of
Hilbert-Schmidt operators:
$$ \pi_2(Tu)=\pi_2(u^*T^*)\leq \|T\|\pi_2(u^*)=\|T\|\|u^*\|=1;
$$
the last equality following from the hypothesis that $X$ has the
$k$-dimensional \twosump\ and the fact that, by (a), $X$ is reflexive.

For (c), let $T:X\to \ell_2^n$ be such that
$\|T\|\|T^{-1}\|=d(X,\ell_2^n)$; then $\sqrt{n}=\pi_2(T^{-1}T)\leq
\|T^{-1}\|\pi_2(T)=\|T^{-1}\|\|T\|=d(X,\ell_n^2)  \leq\sqrt{n}$.
Therefore, $d(X,\ell_2^n)=\sqrt{n}$.\hfill\qed

\bs {\bf Remark.} To make the proof of (b) as simple as possible, we
used (a) to reduce to the case of reflexive spaces.  Actually, it is
well-known that if $\pi_2(T)\le C \nm{T}$ for every operator from $X$
into a Hilbert space $H$, then $X^{**}$ has the same property.  Indeed,
it is easy to see that it is enough to consider finite rank operators
from $X^{**}$ into $H$ and then use local reflexivity (see [LT, p.~33]) and a
weak$^*$ approximation argument.

\proclaim {Example 2.2}  $\ell_{\infty}^2$ has the \twosump.\endproclaim

{\bf Proof:}  Let $u:\ell_\infty^2\to\ell_2^2$, $\|u\|=1$. We can assume that
$$
u=\pmatrix a&b\cr 0&d\cr\endpmatrix
$$
and that $(|a|+|b|)^2+|d|^2\leq1$; or equivalently,
$|a|^2+|b|^2+|d|^2+2|ab| \leq 1$.

For $x=(c_1,c_2)\in\ell_\infty^2$ we have that
$$\leqalignno{
\|ux\|^2=&|ac_1|^2+(|b|^2+|d|^2)|c_2|^2 + 2\Re(ac_1\overline{bc_2})\cr
\le&|ac_1|^2+(|b|^2+|d|^2)|c_2|^2 +|a||b|(|c_1|^2+|c_2|^2).&{(2.1)}}
$$
Set $\lambda=|a|^2+|ab|$,  so that $1-\lambda\geq (|b|^2+|d|^2)+|ab|$.
Thus (2.1) gives
$$
\|ux\|^2\leq |c_1|^2\lambda+|c_2|^2(1-\lambda).\leqno{(2.2)}
$$
Then since (2.2) is in the form of (1.3) with constant 1, we get that
$\pi_2(u)\leq1$.\hfill\qed

\bs

{\bf Remark.}  At least in the complex case, Proposition 2.2 follows from
the fact that $\a(\ell_{\infty}^2)=1$, but we thought it desirable to
give a direct proof.  Another proof is given in [BM].

\ms
Proposition 2.6 provides a useful criterion for determining whether a
space has the {\twosump}.  All of the intuition behind Proposition 2.6
is already contained in Example 2.3:

\proclaim {Example 2.3} Real $\ell_\infty^3$ and complex
$\ell_\infty^4$ do not have the $2$-dimensional \twosump.
\endproclaim

{\bf Proof:} In real $\ell_\infty^3$, let $x_1=(1,0,{1\over\sqrt{2}})$;
$x_2=(0,1,{1\over\sqrt{2}})$ and $X=\span\{x_1,x_2\}$. We denote $X$ by
$X_\infty$ when considered as a subspace of $L_\infty^3$ and by $X_2$
when considered as a subspace of $L_2^3$. Also denote by
$I_{\infty,2}^X$ the restriction to $X$ of the identity $I_{\infty,2}$
from $L_\infty^3$  to $L_2^3$ (we use the standard convention
$L_p^n=L_p^n(\mu)$ where $\mu$ is the probability space assigning mass
${1\over n}$ to every point).

We claim that $\|I_{\infty,2}^X\|<1$ and that $\pi_2(I_{\infty,2}^X)=1$.

For every $\|x\|_\infty=1$,  $\|I_{\infty,2}x\|_2\leq1$ and we have
equality if and only if $|x|$ is {\it flat}; i.e., $x$ is an extreme
point of the unit ball of $L_\infty^3$. Then we verify that
$\|I_{\infty,2}^X\|<1$ by checking that $X$ does not contain any one of
those vectors. For the second one, define $v:\ell_2^2\to X_\infty$ by $v
e_i=x_i$ for $i=1,2$. Then notice that
$\|v\|^2=\|\,|x_1|^2+|x_2|^2\|_\infty=1$, where $|x_1|^2+|x_2|^2$ is
taken coordinatewise in $L_\infty^3$, and  $\pi_2(I_{\infty,2}^X)^2\geq
\pi_2(I_{\infty,2}^X v)^2=\|x_1\|^2_2+\|x_2\|^2_2=1$. The equality
follows, since  $\pi_2(I_{\infty,2}^X) \le \pi_2(I_{\infty,2}) = 1$. We
have thus proved that $X$ does not have the $2$-dimensional {\twosump}.

To conclude, define $u:L_\infty^3\to X_2$ by $u=P I_{\infty,2}$, where
$P$ is the orthogonal projection from $L_2^3$ onto $X_2$.  We claim that
$\|u\|<1$ and that $\pi_2(u)=1$.

If $\|x\|_\infty=1$, then $\|I_{\infty,2}x\|_2=1$ iff $x$ is flat, and
$\|Px\|_2=\|x\|_2$ iff $x\in X$. Since these conditions are mutually
exclusive we conclude that $\|u\|<1$. But $1 =
\|P\|^2\pi_2(I_{\infty,2})^2 \ge \pi_2(u)^2 \geq
\pi_2(uv)^2=\|Px_1\|_2^2+\|Px_2\|_2^2=1$.

The proof for complex $\ell_\infty^4$ is similar: Let
$x_1=(1,0,{1\over\sqrt{2}},{1\over\sqrt{2}})$,
$x_2=(1,0,{i\over\sqrt{2}},{1\over\sqrt{2}})$ and $X=\span\{x_1,x_2\}$.
It is easily checked that $X$ does not contain any flat vectors and
that $|x_1|^2+|x_2|^2\equiv 1$ coordinatewise.\hfill\qed

\medbreak
{\bf Remark.}  We shall
see in section 4 that complex $\ell_\infty^3$ has the $2$-summing
property.
\medbreak

\proclaim {Proposition 2.4}
Let $X$ be an $n$-dimensional subspace of $C(K)$, $K$ compact;
$u:X\to\ell_2^k$ a map satisfying $\pi_2(u)=1$ and $v:\ell_2^k\to X$
satisfying $\|v\|=1$ and $\pi_2(u)=\pi_2(uv)$. Pietsch's factorization
theorem gives the following diagram for some probability $\mu$ on $K$ and
some norm one operator $\a : X_2 \to \ell_2^k$:
$$\matrix
         &                 & C(K)     &{\buildrel I_{\infty,2} \over
                                           \longlongrightarrow}
&L_2(K,\mu)  \cr
         &   &i\uparrow\phantom{i} & &\phantom{P}\downarrow  P  \cr
\ell_2^k & {\buildrel v \over \longrightarrow} & X_\infty &
                           {\buildrel I_{\infty,2}^X \over
                                           \longlongrightarrow}
&X_2        \cr
         &  &\phantom{uuuuu}u\searrow  &
&\swarrow\alpha\phantom{uuuuu} \cr
         &                 &          & \ell_2^k               &
\cr
\endmatrix$$
Let $Y=v(\ell_2^k)$. Then $\alpha$ is an isometry on $Y_2$.
\endproclaim

\medbreak

{\bf Proof:}  We have
$$
1=\nm{uv}^2 =
\sum_{j=1}^k \nm{uve_j}^2 =
\sum_{j=1}^k \nm{\a I_{\infty,2} ve_j}^2 \le
\sum_{j=1}^k \nm{ I_{\infty,2} ve_j}^2 \nm{\a}^2 \le
\pi_2(I_{\infty,2})^2 =1,
$$
so \ $\nm{\a I_{\infty,2} ve_j} = \nm{ I_{\infty,2} ve_j}$ \
for each \ $1\le j \le k$.
Recalling the elementary fact that if $S$ is an operator
between Hilbert spaces, then \ $\{x : \ \nm{Sx} = \nm{S} \nm{x}
\ \}$ \ is a linear subspace of the domain of $S$, we
conclude that \   $\alpha$ is an isometry on $Y_2$. \hfill\qed




\medbreak {\bf Remark.} Let $X$ be an $n$-dimensional space;
$v:\ell_2^n\to X$ the maximum volume ellipsoid map and set
$u={1\over\sqrt{n}}v^{-1}$. It is well known that $\pi_2(u)=1$, and
since $\pi_2(uv)={1\over\sqrt{n}}\pi_2(I)=1$ we conclude that
$\alpha:X_2\to\ell_2^n$
 is an isometry. Moreover, if $X$ is of maximal distance, $u$ is the
minimal volume ellipsoid map (see Remark 1.2).

\medbreak
\proclaim {Corollary 2.5}  Let $X$ be an $n$-dimensional subspace of $C(K)$,
$K$ compact, and $u:X\to\ell_2^n$ be an onto map satisfying $\pi_2(u)=1$.
Suppose that for every orthogonal projection $P$ from $\ell_2^n$ onto a proper
subspace we have $\pi_2(Pu)<1$. Then $\alpha:X_2\to  \ell_2^n$ is an isometry.
($\alpha$ is the map appearing in Proposition 2.4).
\endproclaim

{\bf Proof:} Let $v:\ell_2^n\to X_\infty$ be such that $\|v\|=1$ and
$\pi_2(u)=\pi_2(uv)$ and let $P$ be the orthogonal projection from $\ell_2^n$
onto $u(v(\ell_2^n))=\alpha(Y_2)$. We clearly have that
$\pi_2(Pu)\geq\pi_2(Puv)=\pi_2(uv)=1$. Hence,
the range of $P$ cannot be a
proper subspace of $X_2$ and therefore
$\alpha$ is an isometry on $X_2$.\hfill\qed
\medbreak
In  the next proposition
interpret  ${0\over0}$ as $0$.
 \medbreak

\proclaim {Proposition 2.6} Let $X$ be an $n$-dimensional
subspace of $C(K)$, $K$ compact, and $k \le n$. Then
$$\eqalign{
\sup \biggl\{ & {\pi_2(u)\over\|u\|} : \quad
{u:X_\infty\to\ell_2^k}  \biggr\} =
\cr &
  \sup\biggl\{ { 1
 \over \|P I_{\infty,2}^X\|} \,:\, \mu\in{\cal P}(K), \hbox{
$P^2=P$, $\nm{P}=1$,  rank $P  \le k$,}
\hbox{  $\pi_2(PI_{\infty,2}^X)=1$} \biggr\},  }
$$
where ${\cal P}(K)$ consists of all the probability measures on
$K$, and $I_{\infty,2}$ is the canonical identity from $C(K)$ to
$L_2(K,\mu)$.
\endproclaim

\medbreak

{\bf Proof:}  It is clear that the left hand side dominates the right one. To
prove the other inequality let $u:X_\infty\to\ell_2^k$ be such that
$\pi_2(u)=1$. Then find $v:\ell_2^k\to X_\infty$ such that $\|v\|=1$ and
$\pi_2(uv)=1$. Let $Q$ be the orthogonal projection from $\ell_2^k$ onto
$uv(\ell_2^k)=\alpha(Y_2)$ (with the notation of Proposition 2.4), and $P$ be
the orthogonal projection from $X_2$ onto $Y_2$. Notice that $Qu=\alpha
PI_{\infty,2}^X.$ Since $\alpha$ is an isometry on $Y_2$ we have that
$\pi_2(PI_{\infty,2}^X)=\pi_2(Qu)=1$ and that
$\|PI_{\infty,2}^X\|=\|Qu\|\leq\|u\|$. Therefore
$${1\over\|PI_{\infty,2}^X\|}\geq{\pi_2(u)
\over\|u\|}.$$
\hfill\qed \medbreak

Proposition 2.6 has a nice form when $X$ is $2$-dimensional because then we do
not need to take the orthogonal projection on $X_2$. Indeed, if $P$ has rank
one then it is clear that $\pi_2(PI_{\infty,2}^X)=\|PI_{\infty,2}^X\|$. If
$\mu$ is a probability measure on $K$ with support $K_0$, then
$\|I_{\infty,2}^X\|<1$ iff Ball$(X)$ does not contain any ``flat'' vector on
$K_0$; i.e., whenever $x\in X$ and $\|x\|=1$, then we have that
$|x|_{|K_0}\not\equiv 1$. On the other hand,  $\pi_2(I_{\infty,2}^X)=1$ iff
there exist vectors $x_1, x_2$ in $X$ such that $|x_1|^2+|x_2|^2\leq1$ on $K$
and $|x_1|^2+|x_2|^2\equiv1$ on $K_0$. To see why the second statement is
true, find $v:\ell_2^2\to X_\infty$ satisfying $\|v\|=1$ and
$\pi_2(I_{\infty,2}^X v)=1$. Then let $x_i=ve_i$ for $i=1,2$ and the result is
easily checked for these vectors.  Since every closed subset of a
compact metric space is the support of some probability measure, this
discussion proves:

\medbreak
\proclaim {Corollary 2.7.a}   Let $X$ be a $2$-dimensional subspace of
$C(K)$, $K$ compact metric. Then $X$ does not have the $2$-summing property
if and only if there exist vectors $x_1, x_2$ in $X$ and a closed set
$K_0\subset K$ with $1_{K_0} \le  |x_1|^2+|x_2|^2\leq 1$
 such that for every $x\in X$ with
$\|x\|_\infty=1$, we have that $|x|_{|K_0}\not\equiv 1$.
\endproclaim

\medbreak

If $X\subset C(K)$ contains a vector $|x|\equiv 1$, then for every probability
measure $\mu$ on $K$ we have $\|I_{\infty,2}x\|=\|x\|$; hence,
$\|I_{\infty,2}^X\|=1$ and we have

\proclaim {Corollary 2.8} Let $X$ be a $2$-dimensional subspace of $C(K)$, $K$
compact. If $X$ contains a vector $x$, $|x|\equiv 1$ on $K$ then $X$ has the
$2$-summing property.
\endproclaim

This applies immediately to $\ell_\infty^2$ (both real and complex) and also
to $\ell_1^2$ (again real and complex) if  embedded in a canonical way.  It
also implies that there are a continuum of pairwise nonisometric two
dimensional complex spaces which have the {\twosump}.
We shall see later that the real \twosump\  is quite different
from the complex version. For the moment,  take
$X=\span\{(1,0,{1\over\sqrt{2}}),(0,1,{1\over\sqrt{2}})\}$ inside
$\ell_\infty^3$. We proved in Example 2.3 that real $X$ does not have the
\twosump. However, complex $X$ has it. To see this, notice that
$(1,0,{1\over\sqrt{2}})+i(0,1,{1\over\sqrt{2}})$ is ``flat'' and hence
Corollary 2.8 implies the result. The difference can be explained by
saying that it is easier to get ``flat'' vectors in the complex
setting (see Proposition 4.4).

Let us return to the discussion following Proposition 2.6. Suppose that
$X$ is an $n$-dimensional subspace of $C(K)$, $\mu$ is a probability
measure on $K$ with support $K_0\subset K$, and $P$ is an orthogonal
projection from $X_2 \subset L_2(\mu)$ onto a subspace $Y_2$. Notice that
\ $\pi_2(PI_{\infty,2}^X)=1$ \ if and only if there exist vectors
$x_1,x_2,\dots,x_n$ in $Y$ with  $1_{K_0} \le \sum_{j=1}^n |x_j|^2 \le
1$ and each vector $I_{\infty,2}x_j$ is in $Y_2$.
(Keep in mind that   $Y_2$ is relatively $L_2(\mu)$--closed in $X_2$,
hence if $y\in Y$, $z\in X$, and $1_{K_0}y=1_{K_0}z$, then also $z$ is in
$Y$.)
Similarly,
\ $\nm{PI_{\infty,2}^X}=1$ \ if and only if there exists a {\sl
single\/} vector $x$ in $Y$  with \  $1_{K_0}\le |x|\le 1$. Thus we get a
version of Corollary 2.7.a for all finite dimensional spaces:

\proclaim {Corollary 2.7.b}   Let $X$ be a finite dimensional subspace
of $C(K)$, $K$ compact metric. Then $X$ does not have the $2$-summing
property if and only if there exist vectors $x_1, x_2, \dots, x_n$ in $X$
and a closed set  $K_0\subset K$ with
$1_{K_0} \le \sum_{j=1}^n |x_j|^2 \le 1$ such that for every $x\in X$ with
$1_{K_0} \le |x|\le 1$, we have that $1_{K_0}x$ is not in
$\span\{ 1_{K_0}x_1,1_{K_0}x_2,\dots,1_{K_0}x_n\}$.
\endproclaim

\medbreak

In the complex case, the \twosump\ is not hereditary, since complex
$\ell_1^3$ has the \twosump\ but $\ell_1^2$ is the only two
dimensional subspace of it which has the \twosump\ (see Theorem 4.2
and Proposition 5.7.) Nevertheless:

\proclaim{Proposition 2.9} Let $X$ be a subspace of $\ell_\infty^N$
which has the \twosump. Then every subspace of $X$ has the \twosump.
\endproclaim

{\bf Proof:} Assume that $X$ has a subspace  which fails the
\twosump. Write $K=\{1,2,\dots,N\}$. Since $\ell_\infty^N=C(K)$, we can apply
Corollary 2.7.b. There exists a subset  $K_0\subset K$ for which we can find
vectors $x_1, x_2, \dots,x_n$ in $X$  with  $1_{K_0} \le
\sum_{j=1}^n |x_j|^2 \le 1$ such that no norm one vector in
$Y\equiv\span\{\seq{x}{n}\}$ is
unimodular on $K_0$. We can also assume that $K_0$ is maximal
with respect to this property; in particular,
$\sum_{j=1}^n |x_j|^2$ is strictly less than one off $K_0$ and hence
$\sum_{j=1}^n 1_{\sim K_0}|x_j|^2 < 1-\e$ for some $\e>0$.

On the other hand,  since $X$ has the
\twosump, there exists a vector $y\in Y$ which is unimodular on $K_0$
(and whose restriction to $K_0$ agrees with the restriction to $K_0$ of
some unit vector in $X$).  Evidently $\nm{y}>1$.  Thus there exists
$1>\tau>0$ so that
$z=z_\tau\equiv  |\tau^{1\over2}y|^2 + \sum_{j=1}^n
|(1-\tau)^{1\over2}x_j|^2$ satisfies
$\nm{1_{\sim K_0}z}_\infty = 1$.  But then $z\le 1$, a
square function of a system from $Y$,  is unimodular on a set which
properly contains $K_0$; this contradicts the maximality of $K_0$.
\hfill\qed

\ms

For any real Banach space $F$, let $F_\C$ denote
the linear space $F\oplus F$ with complex structure
defined by means of the formula
$(a+bi)(f_1\oplus f_2)=(af_1-bf_2)\oplus(af_2+bf_1)$.
There is a natural topology on $F_\C$, namely $F_\C$ is
homeomorphic with the direct sum of two real
Banach spaces $F\oplus F$. In two special cases
we shall define a norm on $F_\C$ which will make it a
complex Banach space. First, if $E$ is a linear subspace
of some  real $\ell_\infty^k$, we endow $E_\C$ with the norm
induced from complex $\ell_\infty^k$ by means of the obvious
embedding.  Secondly, if $E=H$ is a Hilbert space, then
$H_\C$ is normed by means of the formula $\nm{h_1\oplus
h_2}=(\nm{h_1}^2+\nm{h_2}^2)^{1/2}$. These two
definitions are consistent, because $H$ is isometric to a
subspace of  real $\ell_\infty^k$  only if $\dim_{\R} H\le1$.
Now, if $T:F\rightarrow G$ is a linear operator, we define
$T_\C:F_\C\rightarrow G_\C$ by the formula
$T_\C(f_1\oplus f_2)=(T f_1\oplus T f_2)$.
\medbreak

\proclaim {Proposition 2.10} Let $E$ be a subspace of
real $\ell_\infty^k$ and $S:E\rightarrow \ell_2^2$ be a real-linear
mapping from $E$ into a $2$-dimensional real Hilbert space. Then
$$
\pi_2(S_\C)=\pi_2(S)=\nm{S_\C}.
$$
\endproclaim

{\bf Proof:}   From the discussion in Section 1 we see that there are
vectors $x,y$ in $E$ such that $|x|^2 + |y|^2\le 1$ (interpreted
coordinatewise) and $$ \pi_2(S)=\bigl(\nm{Sx}^2+\nm{Sy}^2\bigr)^{1/2}. $$

Observe that, by our definition, $\nm{x\oplus y}_{E_\C}=
\nm{x+iy}_{\l_\infty^k(\C)}\le1$. It follows that  $$
\pi_2(S_\C)\ge\nm{S_\C}\ge\nm{S_\C(x\oplus y)}_{H_\C}=
\bigl(\nm{Sx}^2+\nm{Sy}^2\bigr)^{1/2}=\pi_2(S), $$ hence it remains to
verify that $\pi_2(S_\C)\le\pi_2(S)$.

Take $u,w$ in $E_\C$ with $|u|^2+|w|^2\le 1$ and
$\pi_2(S_\C)^2=\Vert{S_\C(u)}\Vert^2+\Vert{S_\C(w)}\Vert^2$. Interpreting
real and imaginary parts of vectors in $\ell_\infty^k$ coordinatewise, we
see that $\Re{u},\Im{u}$ are in $E$ and similarly for $w$.  Moreover,  $$
\Vert{S_\C(u)}\Vert^2+\Vert{S_\C(w)}\Vert^2=
\nm{S\Re{u}}^2+\nm{S\Im{u}}^2+\nm{S\Re{w}}^2+\nm{S\Im{w}}^2. $$ This last
quantity is at most $\pi_2(S)$ since
$|\Re{u}|^2+|\Im{u}|^2+|\Re{w}|^2+|\Im{w}|^2 = |u|^2+|w|^2\le 1$. \hfill\qed
\medbreak
It is easy to determine when a complex-linear operator is the
complexification of a real-linear operator:
\ms

\proclaim {Proposition 2.11} Let $E$ be a real Banach space and let
$G$ be a complex Hilbert space. Let $T:E_\C\rightarrow G$ be a
complex-linear continuous operator. The following conditions
are equivalent:
\item{(i)} There is a real Hilbert space $H$ and continuous
linear operators $S:E\rightarrow H$, $U:H_\C\rightarrow G$, such
that $T=U\circ S_\C$ and $U$ is an isometric embedding.
\item{(ii)} For each $e,e'\in E$ one has
$\nm{T(e+ie')}=\nm{T(e-ie')}$.
\item{(iii)} For each $e,e'\in E$ one has $\Im(Te,Te')=0$.
\item{(iv)} There is a subset $E_0$ of $E$ such that the linear
span of $E_0$ is dense  in $E$ and for each $e$, $e'$ in $E_0$,
 $\Im(Te,Te')=0$.
\endproclaim

{\bf Proof:} (i) implies (ii), because $\nm{T(e\pm ie')}=
\nm{S_\C(e\pm ie')}=(\nm{Se}^2+\nm{Se'}^2)^{1/2}$.
To see that (ii) implies (i) we let $H$ denote the closure
of $T(E)$ in $G$. Observe that, if $x=e\oplus e'\in E_\C$,
then using the parallelogram identity we obtain
$$
\nm{Te+iTe'}^2={1\over2}(\nm{Te+iTe'}^2+
\nm{Te-iTe'}^2)=\nm{Te}^2+\nm{Te'}^2.
$$
Hence $H\oplus iH$ is linearly isometric to $H_\C$
and, if $U$ denotes the natural embedding and $S=T|_E$,
then we have $T=U\circ S_\C$.

It is clear that (iii) and (iv) are equivalent. On the
other hand, the identity
$$
\eqalign{
\nm{T(e+ie')}^2-\nm{T(e-ie')}^2&=
(Te+iTe',Te+iTe')-(Te-iTe',Te-iTe')\cr
&=2i((Te',Te)-(Te,Te'))=4\Im(Te,Te'),}
$$
makes it obvious that (ii) and (iii) are equivalent. \hfill\qed

\medskip

\bs

\head 3. The real case.\endhead
Throughout this section we deal with spaces over the
reals.  Example 2.3 and Proposition 2.1 imply that there are many norm
one operators from real $\ell_1^3$ into $\ell_2^2$ whose \twosum\ norm is
larger than one.

\proclaim {Proposition 3.1}  Let $u$ be an operator from real
$\ell_1^3$ to $\ell_2^2$  such that  $ue_1, ue_2, ue_3$  have norm one and
every two of them are linearly independent. Then $\pi_2(u)>1$.
\endproclaim

{\bf Proof:}  Let $u:\ell_1^3\to\ell_2^2$ satisfy the assumptions of
Proposition 3.1,  and set
$x_i=ue_i$, $i=1,2,3$.  Notice that $\ell_1^3$ embeds isometrically
via the natural evaluation mapping into $C(K)$, where
$K=\{(1,1,1),(1,1,-1),(1,-1,1),(-1,1,1)\})$, regarded as a subset of
$\ell_\infty^3 = (\ell_1^3)^*$.
 Assume that $\pi_2(u)=1$ and consider the Pietsch factorization
diagram:

$$\matrix
         &                 & C(K)     &{\buildrel I_{\infty,2} \over
                                           \longlongrightarrow}
&L_2(K,\lam)  \cr
         &   &i\uparrow\phantom{i} & &\phantom{P}\downarrow  P  \cr
         &                         & X_\infty &
                           {\buildrel I_{\infty,2}^X \over
                                           \longlongrightarrow}
&X_2        \cr
         &  &\phantom{uuuuu}u\searrow  &
&\swarrow\alpha\phantom{uuuuu} \cr
         &                 &          & \ell_2^2               &
\cr
\endmatrix$$

So $\nm{\a}=\pi_2(u)=1$.  We claim that $\dim X_2=3$. (This is not
automatic since the support $K_0$ of $\lam$ may not be all of
$K$.) Indeed, any two of $\{f_1,f_2,f_3\}$ are linearly independent,
since this is true of $\{\a f_1,\a f_2,\a f_3\}$, so the cardinality
of $K_0$ is at least three.  But
$\{{f_i}_{|L}\}_{i=1}^3$ is linearly independent if $L$ is any three
(or four) element subset of $K$.

To complete the proof, just recall that no norm one linear operator
$\b$ from $\ell_2^3$ into $\ell_2^2$ achieves its norm at three
linearly independent vectors since, e. g., the set \
$\{x\in\ell_2^3 :  \nm{\b x}=\nm{\b}\nm{x} \}$ \  is a subspace of
$\ell_2^3$.\hfill\qed

\ms

Having treated the ``worst" case of $\ell_1^n$, it is easy to formulate a
version of Proposition 3.1 for general spaces.

\ms
\proclaim {Corollary 3.2} Let $X$ be a real space. \item{(a)} If $T$ is  a
norm one operator from $X$ into $\ell_2^2$ such that for some $z_1,z_2,z_3$
in  $\Ball(X)$, $Tz_1, Tz_2, Tz_3$ have norm one and every two of them  are
linearly  independent, then $\pi_2(T)>1$. \item{(b)}  If there exist
$u:\ell_2^2\to X$  and $x_1,x_2,x_3\in\Ball(\ell_2^2)$ such that
$\|u\|=1=\|ux_i\|=1$ for  $i=1,2,3$ and every two of $x_1,x_2,x_3$ are
linearly independent, then $X$ does not have the $2$-dimensional \twosump.
\endproclaim

{\bf Proof:} For (a), define $w:\ell_1^3\to X$ by $we_i=z_i, i\leq3$. Then
 $Tw$ satisfies the
hypothesis of Proposition 3.1. Therefore $\pi_2(Tw)>1$; and since $\|w\|=1$
we have that
$\pi_2(T)>1$.

For (b) it is enough to prove that  $X^*$ does not
have the $2$-dimensional \twosump.
 If $x_i^*\in \Ball(X^*)$, $i=1,2,3$
satisfy $\scpr{x^*_i,ux_i}=1$,
then
$1\ge \nm{u^*x_i^*}\ge \scpr{u^*x_i^*,x_i}=\scpr{x^*_i,ux_i}=1$; therefore,
$u^*x_i^*=x_i$ for $i=1,2,3$ and the previous part
gives us that $\pi_2(u^*)>1$. \hfill \qed

\ms

We are now ready for the main result of this section:

\proclaim {Theorem  3.3} If $X$ is a real space of dimension at
least three, then $X$ does not have the \twosump.  Consequently,
the only real spaces which have the \twosump \ are $\R$ and
$\ell_{\infty}^2$.
\endproclaim

{\bf Proof:}  The ``consequently" follows from the first statement and
Proposition 2.1(c) because in the real case $\ell_{\infty}^2$ is the only
$2$-dimensional maximal distance space.  This is an unpublished result of
Davis and the second author; for a proof see Lewis' paper [L].

So assume that $X$ is $\R^n$, $n\ge 3$,  under some norm and has the
\twosump.  We can also assume that the usual Euclidean ball $\cal E$ is the
ellipsoid of minimal volume containing $\Ball(X)$, and we use $|\cdot|$ to
denote the Euclidean norm.

By John's theorem, there exist $\mu_1,\cdots,\mu_m > 0$ such that
$\sum_{i=1}^m\mu_i=n$  and
$y_1,\cdots,y_m\in X$ outside contact points (i.e., $\|y_i\|=|y_i|=1$
for $i=1,2,\dots,m$) such
that every $x\in X$ satisfies $x=\sum_{i=1}^m\mu_i\scpr{x,y_i}y_i$.  Recall
that ${\cal E}/\sqrt{n}\subset\Ball(X)$, in fact, by Lemma 1.1,
${\cal E}/\sqrt{n}$ is the ellipsoid
of maximal volume contained in $\Ball(X)$.

If $x$ is an inside contact point (i.e., $\|x\|=1$ and $|x|=1/\sqrt{n}$),
Milman and  Wolfson [MW] proved that
$$
|\scpr{x,y_i}| = {1\o n} \quad\hbox{  for every }
i=1,2,\dots,m.\leqno{(3.1)}
$$
To see this, observe that $\{z\in X : \scpr{z,x} = {1/n} \}$ supports ${\cal
E}/\sqrt{n}$ at the inside contact point $x$, hence--draw a picture--the norm
of $\scpr{\cdot,x}$ in $X^*$ is $1/n$. Thus
$$
{1\o n} = \scpr{x,x} = \sum_{i=1}^m \mu_i\scpr{x,y_i}^2
\le \sum_{i=1}^m \mu_i\nm{\scpr{\cdot,x}}_{X^*}^2\nm{y_i}^2={1\o n}.
$$
This implies that $|\scpr{x,y_i}| = {1\o n}$ for $i=1,2,\dots,m$ and proves
(3.1).

In other words, the norm one (in $X^*$) functional  $\scpr{\cdot,nx}$ norms
all of the $y_i$'s as well as $x$.   This implies  that both
conv$\{y_i:\scpr{y_i,nx}=1\}$ and conv$\{y_i:\scpr{y_i,nx}=-1\}$ are subsets
of the unit sphere of $X$.  Now we know that $X$ has maximal distance, hence
at least one inside contact point exists, whence $\Ball(X)$ has at least two
``flat" faces.

The next step is to observe that we can find $n$ linearly independent outside
contact points  $y_1,\cdots,y_n$ and $n$ linearly independent inside contact
points  $x_1,\cdots,x_n$ satisfying (3.1).  This will give us enough faces on
$\Ball(X)$ so that a  $2$-dimensional section will be a hexagon and we can
apply Corollary 3.2 (b).  Now the John representation of the identity gives
the existence of the outside contact points, and since  ${\cal E}/\sqrt{n}$
is the ellipsoid of maximal volume contained in $\Ball(X)$, another
application of John's theorem gives the inside contact points.

So let $y_1,\cdots,y_n$ be linearly independent outside contact points, and
$x_1,\cdots,x_n$ linearly independent inside contact points satisfying (3.1),
such that for the first three of them we have, $$
\eqalign{\scpr{nx_1,y_1}=1&\quad\scpr{nx_1,y_2}=1\quad \ \
\scpr{nx_1,y_3}=1\cr
           \scpr{nx_2,y_1}=1&\quad\scpr{nx_2,y_2}=1\quad \ \
\scpr{nx_2,y_3}=-1\cr
           \scpr{nx_3,y_1}=1&\quad\scpr{nx_3,y_2}=-1\quad
\scpr{nx_3,y_3}=1.\cr}  $$ (We are allowed to change signs and renumber the
contact points). Let $v$ denote the linear map from $\span \{y_1,y_2,y_3\}$
into $l_1^3$ which takes $y_i$ to $e_i$.  Then $\nm{vy}_1 = \nm{y} $ if
$y=\sum_1^3 \a_i y_i$ and either $\a_1\a_2 \ge 0$ or $\a_1\a_3 \ge 0$, hence
the restriction of $v$ to  $F = \span\{{{y_1+y_2}\over 2}, {{y_1+y_3}\over
2}\}$ is an isometry.  But then $v[\Ball(F)]$ is a regular hexagon, which
implies that the maximum volume ellipsoid for $\Ball(F)$ touches the unit
sphere of $F$ at six points.  We finish by applying   Corollary 3.2
(b).\hfill\qed

\medskip Lewis [L] proved that every real maximal distance space of dimension
at least two contains a subspace isometrically isomorphic to $\ell_1^2$; that
is, a subspace whose unit ball is a parallelogram.  In view of Remark 1.2,
the proof of Theorem 3.3 yields:

\medskip

\proclaim {Corollary 3.4} If $X$ is a real maximal distance space of
dimension at least three, then $X$ has a subspace whose unit ball is a
regular hexagon.
\endproclaim

\bs

\head 4. The complex case.\endhead
As was mentioned in the introduction, the study of
the \twosump \  in the complex case is {\sl a priori\/} more complicated than
in the real case even for two dimensional spaces simply because in the real
case there is only one $2$-dimensional maximal distance space, while in the
complex space there are at least two, $\ell_{\infty}^2$ and its dual.  In
fact, it is not difficult to construct other complex $2$-dimensional maximal
distance spaces without using Corollary 2.8.  One way is to use John's
representation theorem; this was done independently by Gowers and
Tomczak-Jaegermann [both unpublished] a couple of years ago in order to
construct real $4$-dimensional maximal distance spaces whose unit ball is not
strictly convex; this approach also yields maximal distance $2$-dimensional
complex spaces different from $\ell_{\infty}^2$ and $\ell_{1}^2$.  However, a
simpler way of seeing that there are many  maximal distance $2$-dimensional
complex spaces is via Proposition 4.1:

\proclaim {Proposition 4.1} Suppose that $X$ is an $n$-dimensional complex
space which, as a real space, contains a real $n$-dimensional subspace which
has maximal distance to real $\ell_2^n$.  Then $X$ has maximal distance to
complex $\ell_2^n$.
\endproclaim

{\bf Proof:}  Let $Y$ be a real $n$-dimensional subspace of $X$ which  has
maximal distance to real $\ell_2^n$, and suppose that $T$ is any complex
linear isomorphism from $X$ to complex $\ell_2^n$. But considered as a real
space, complex $\ell_2^n$ is just real $\ell_2^{2n}$, and so the restriction
of $T$ to $Y$ is a real linear isomorphism from $Y$ to a real $n$-dimensional
Hilbert space, hence by the assumption on $Y$,  $$ \nm{T}\nm{T^{-1}}\ge
\nm{T_{|Y}}\nm{(T_{|Y})^{-1}}\ge \sqrt{n}, $$ so the desired conclusion
follows.\hfill\qed

\ms

Proposition 4.1 makes it easy to construct in an elementary manner
$2$-dimensional complex maximal distance spaces which are not isometric to
either $\ell_{\infty}^2$ or $\ell_{1}^2$.  For example, in complex
$\ell_{\infty}^3$, let $x=(1,0,a)$ and $y=(0,1,b)$ with $|a+b|\vee |a-b|\le
1$ but $|a|+|b| > 1$, and set $X=\span\{x,y\}$.  In Proposition 4.4  we
prove that {\sl every\/} $2$-dimensional subspace of complex
$\ell_{\infty}^3$ is a maximal distance space, by proving that they all
have the \twosump.

However, in the complex setting, the \twosump \  is not restricted to two
dimensional spaces.  In fact, we do not know a good bound for the dimension of
complex spaces which have the \twosump, although we suspect that dimension
three is the limit.  In dimension three itself, we know of only two examples,
$\ell_1^3$ and its dual.

We can prove Theorem 4.2 using Proposition 2.6 and The proof of
Proposition 4.4 (see Remark 4.6). The proof we give is of independent
interest.

\proclaim {Theorem 4.2} Complex  $\ell_1^3$  has the \twosump.  Hence also
$\ell_{\infty}^3$  has the \twosump.
\endproclaim

{\bf Proof:}  The proof reduces the theorem to the following calculus lemma,
which we prove after giving the reduction:

\ms

\proclaim {Lemma 4.3}  Fix arbitrary complex numbers $\lambda_1$,
$\lambda_2$, $\lambda_3$ and define a function $f$ on the bidisk by $$
f(z_1,z_2) = |1 + \lambda_1 z_1 + \lambda_2 z_2 + \lambda_3 z_1 \bar{z_2}| +
|\lambda_3| \sqrt{1-|z_1|^2}\sqrt{1-|z_2|^2}.   $$ Then the maximum of $f$
is attained at some point on the two dimensional torus; that is, when
$|z_1| = |z_2| = 1$.
\endproclaim

\ms

{\bf Reduction to Lemma 4.3:} Let $K= \{1\}\times \T \times \T$, where $\T$
is the unit circle, and regard $\ell_1^3$ as the subspace of $C(K)$ spanned
by the coordinate projections $f_0, f_1, f_2$.  Let $u$ be a norm one
operator from $\ell_1^3$ into $\ell_2$ (complex scalars).  We want to show
that $\pi_2(u) = 1$.  That is, we want a probability $\mu$ on $K$ so that
for each $x \in \ell_1^3$,
$$
\nm{ux}^2 \le \int |x(k)|^2\, d\mu(k).
$$
Notice that if we add to the $u f_j$'s mutually orthogonal vectors
which are also orthogonal to the range of $u$, the $\pi_2$-norm of the
resulting operator can only increase.  Thus we assume, without loss
of generality, that
$\nm{u f_j} = 1$ for $j=0,1,2$,  and   that
$$
uf_0 = \d_0, \quad
uf_1 = \a_1\d_0 + \b_1 \d_1, \quad  uf_2 = \a_2\d_0 + \b_2 \d_1 + \g_2
\d_2,
$$
where $\{\delta_0,\delta_1,\delta_2\}$ is an orthonormal set of $\ell_2$.
Define a linear functional $F$ on $E=\span \{f_0,f_1,f_2,f_1\bar{f_2}\}$ by
$$
Ff_0 = 1,\quad Ff_1 = \a_1 \quad Ff_2 = \a_2,
\quad F(f_1\bar{f_2}) = \a_1 \bar{\a_2} + \b_1\bar{\b_2}.
$$

{\bf Claim:}  $\nm{F} = 1$ {\sl as a linear functional on }
$(E,\nm{\cdot}_{C(K)})$.

\medskip

Assume the claim.  Then by the Hahn-Banach theorem, $F$ can be extended to a
norm one linear functional, also denoted by $F$, on $C(K)$. Since $\nm{F} =
1=Ff_0$,  $F$ is given by integration against a probability, say, $\mu$.  Now
by the definition of $F$, the mapping
$v :\{0,f_0,f_1,f_2\}\to \ell_2$ defined by
$v 0 = 0$, \ $v f_j = u f_j$  for $j=0,1,2$ is
$L_2(\mu)$-to-$\ell_2$ inner product preserving, hence an $L_2$-isometry,
whence extends to a linear isometry from
$(\span \{f_0,f_1,f_2\},\nm{\cdot}_{L_2(\mu)})$ into $\ell_2$.  This shows that
$\pi_2(u)\le 1$, as desired.

\medskip

{\bf Proof of claim:}  The claim says that for all
$\{\lam_j\}_{j=0}^3 $ in $\C^{4}$,
$$
|\lam_0 + \lambda_1 \a_1 + \lambda_2 \a_2 +
\lambda_3 ( \a_1 \bar{\a_2} + \b_1 \bar{\b_2})| \le
\sup_{(z_1,z_2)\in T\times T}
|\lam_0 + \lambda_1 z_1 + \lambda_2 z_2 + \lambda_3 z_1 \bar{z_2}|.
$$
But the left hand side is dominated by
$$
|\lam_0 + \lambda_1 \a_1 + \lambda_2 \a_2 +
\lambda_3 \a_1 \bar{\a_2} | + |\lambda_3| \sqrt{1-|\a_1|^2}\sqrt{1-|\a_2|^2},
$$
so the claim follows from Lemma 4.3.

\ms

{\bf Proof of Lemma 4.3:}  For fixed $z_2$, you can rotate $z_1$ to see that at
the maximum
$$
f(z_1,z_2) =  |1 +  \lambda_2 z_2 | +
|z_1| |\lambda_1 + \lambda_3 \bar{z_2}| +
|\lambda_3| \sqrt{1-|z_1|^2}\sqrt{1-|z_2|^2}.
$$

As $|z_1|$ varies, $(|z_1|,\sqrt{1-|z_1|^2})$ varies over the unit sphere of
$\ell_2^2$, so at the maximum
$$
f(z_1,z_2) = |1 +  \lambda_2 z_2 | + \sqrt{|\lambda_1 + \lambda_3
\bar{z_2}|^2  + |\lambda_3|^2 (1-|z_2|^2)}.
$$
In particular, $|z_1| = 1$ if and only if $|z_2|=1$.  The last expression can
be rewritten as
$$
|1 +  \lambda_2 z_2 | +  \sqrt{|\lambda_1|^2 + |\lambda_3|^2 +
2\Re(\lambda_1 \bar{\lambda_3} z_2)}.
$$
Choose $\theta$ so that $e^{i\theta} \lambda_1 \bar{\lambda_3}$ is purely
imaginary; so ${\Re }(\lambda_1 \bar{\lambda_3} z_2)$ does not change if you
add a real multiple of $e^{i\theta}$ to $z_2$.  Assume first that
$\lambda_2\neq0$. Then, for any $\epsilon > 0$,
by adding  either $\epsilon e^{i\theta}$ or $- \epsilon e^{i\theta}$ to $z_2$
you increase the first term; this means that at the maximum
$|z_2| = 1$. The case where $\lambda_2=0$ can be handled in a similar way.
\hfill\qed

\ms

{\bf Remarks.} 1. A similar argument reduces the problem of whether
$\ell_1^4$ has the \twosump\  to a calculus problem; however, in the remark
 after Example 2.3 we give a simple direct argument that $\ell_{\infty}^4$
fails the \twosump.

2.  Bagchi and Misra [BM] give a  different reduction of Theorem 4.2 to a
variation of Lemma 4.3.  Their argument may be more appealing to operator
theorists.

\ms

3. The following Proposition is a consequence of Proposition 2.9 and
Theorem 4.2. We indicate a shorter argument.
\ms

\proclaim {Proposition 4.4} Every two dimensional subspace  of
$\ell_\infty^3$ has the \twosump.
\endproclaim

{\bf Proof:} Assume that $X\subset\ell_\infty^3$ is two
dimensional. Applying  $\ell_\infty^3$ isometries we assume that $X$ has a
basis of the form $\{(1,0,a_1), (0,1,a_2)\}$ where $|a_i|\leq1$ for
$i=1,2$. If $|a_1|+|a_2|\leq1$ then $X$ is isometric to $\ell_\infty^2$ and
by Example 2.2 it has the \twosump. So assume that $|a_1|+|a_2|>1$. It is
easy to find $\theta$, $0\leq\theta<2\pi$ such that
$|a_1+e^{i\theta}a_2|=1$. Hence, $(1,0,a_1)+e^{i\theta}(0,1,a_2)$ is
``flat'' and the result follows from Corollary 2.8.\hfill\qed

\medbreak

It is interesting to notice that Proposition 4.4 is equivalent to the following
 calculus formulation.

\proclaim {Proposition 4.5} Given any
complex numbers $c_1,c_2,c_3$ and $d_1,d_2,d_3$  with $|c_j|^2 +
|d_j|^2 \le 1$ for $j=1,2,3$, suppose that  $\a,\b$ satisfy
$$
|\a|^2 |\gamma|^2 + |\b|^2 |\d|^2 \le
\max_{j=1,2,3} |\gamma c_j + \d d_j|^2 \quad \forall \gamma, \d \in
\C. \leqno{(4.1)}
$$
Then $|\a|^2 + |\b|^2 \le 1$.
\endproclaim

To see the equivalence of Propositions 4.4 and 4.5,  let $X$ be a
$2$-dimensional subspace of $\ell_\infty^3$ and  $X {\buildrel u \over
\longrightarrow} \ell_2^2$ a norm one operator.  Choose  $\ell_2^2
{\buildrel v \over \longrightarrow} X$  of norm one so that   $\pi_2(u) =
\pi_2(uv)$.  We can assume that  $uv$ is diagonal; say,  $uv(e_1) = \a e_1$
and $uv(e_2) = \b e_2$.  For $j = 1,2$ set $x_j = v(e_j)$ and write $x_1 =
(c_1,c_2,c_3)$,    $x_2 = (d_1,d_2,d_3)$.  So  $\pi_2(u)^2 = |\a|^2 +
|\b|^2 $, while $\nm{v}^2 ={\displaystyle \max_{j=1,2,3}}|c_j|^2 +
|d_j|^2 = 1$. The implication Proposition 4.5 $\Rightarrow$ Proposition 4.4
follows by noticing that $u$ having norm at most one is equivalent to the
inequality (4.1).  Similar considerations yield the easier reverse
implication.

We do not see a really simple proof of the calculus reformulation of
Proposition 4.4 without using Pietsch's factorization theorem.  However, a
 similar reduction of
the weaker statement that every $2$-dimensional subspace of
$\ell_\infty^3$ has maximal distance to $\ell_2^2$ produces a calculus
statement which is very easy to prove.  Indeed,  given a $2$-dimensional
subspace $X$ of $\ell_\infty^3$, we can choose norm one operators $X
{\buildrel u \over \longrightarrow} \ell_2^2$ and $\ell_2^2 {\buildrel v
\over \longrightarrow} X$ so that $ d uv = I_{\ell_2^2}$, where
$d$ is the distance from $X$ to $\ell_2^2$.  We can choose the
orthonormal basis $e_1,e_2$ so that $1=\nm{v}=\nm{ve_1}$ and define
$x_1,x_2$, the $c_i$'s,and the $d_i$'s as in the discussion above.  Since
$\nm{x_1}=1$, we can assume, without loss of generality, that
$|c_1|=1$.  (4.1) holds with $\a=\b={1\over d}$, and we want to see that
this implies that $d\ge \sqrt 2$; i.e., that $\a^2 + \b^2 \le 1$. So we
only need to get $\g $ and $\d $ of modulus  one to make the right side
of (4.1)  one.  Since $|c_1| = 1$, $d_1=0$, any such choice makes
$|\g  c_1 + \d  d_1| = 1$. Choose $\g $ to make  $\g  c_2 \ge 0$; then
$|\g  c_2 + \d d_2| \le 1$ as long as $\d d_2$ has nonpositive real part,
which happens as long as $\d$ is on a certain closed semicircle.
Similarly,  $|\g  c_3 + \d d_3| \le 1$ as long as $\d $ is on another
closed semicircle.  Since any two closed semicircles of the unit circle
intersect, the desired choice of $\g $ and $\d $ can be made.

\medskip

{\bf Remark 4.6.}  1. We can use Propositions 2.6 and 4.4 to prove that
$\ell_\infty^3$ has the $2$-summing property. Use the notation of
Proposition 2.4 and let $P$ be an orthogonal projection on $L_2^3(\mu)$. If
the rank is one, there is nothing to prove. If the rank is three then we
clearly have that  $\|I_{\infty,2}\|=1$, and if the rank is two, then the
proof of Proposition 4.4 implies the result.

2. It is natural to ask if the only subspaces of complex $L_1$ with the
\twosump\ are  $\ell_1^2$ and $\ell_1^3$. The answer is yes
because a subspace of $L_1$ of maximal distance is already an $\ell_1^k$
space. We prove this in the appendix, Proposition 5.7.

\bs

\head 5. Appendix.\endhead

In this section we present some related results.

\proclaim {Proposition 5.1} Every subspace of complex $\ell_\infty^3$ is
the complexification of a subspace of real $\ell_\infty^3$
\endproclaim

The proof of Proposition 5.1 follows easily from the next two lemmas.
Recall that a vector in $\ell_\infty^k$ is said to be {\sl flat \/} if all
of its coordinates are unimodular.

\proclaim {Lemma 5.2} Suppose that the subspace  $X$ of
complex $\ell_\infty^3$ is not linearly isometric to
$l_\infty^2$. Then $X$ contains two linearly independent flat
vectors, say $f_1$ and $f_2$. Moreover, each flat vector in $X$
is of the form $\lambda f_j$, where $j\in\{1,2\}$ and
$|\lambda|=1$.
\endproclaim

{\bf Proof:}
Applying $\ell^3_\infty$ isometries, we may assume that $X$ is spanned by
two vectors of the form $x=(1,0,a)$ and $y=(0,1,b)$ where $a,b\in\C $
with $\vert a\vert,\vert b\vert \le 1$. Put $w= x-\psi y$. For $w$
to be flat one needs that $|\psi|=1$ and $|a-\psi b|=1$.
Observe that, since $X$ is not linearly isometric to $l_\infty^2$, we
have $|a|+|b|>1$;, in particular  $ab\ne0$. Thus $\psi\in\C$ should
belong to the intersection of the unit circle $\{z\,:\,|z|=1\}$ and the
circle  $\{z\,:\,|z-a/b|=1/|b|\}$, hence there are at most two solutions
for $\psi$. Thus it will suffice to check that the two circles have a
point in common and that they are not tangent at that point. Since the
second circle has a bigger radius, this amounts to verifying the strict
inequalities   $$ \Bigl|{1\over b}\Bigr|-1<\Bigl|{a\over b}-0\Bigr|
<1+\Bigl|{1\over b}\Bigr|.   $$ These inequalities are obvious, because
we have $|a|+|b|>1$, $|a|\le1$ and $|b|>0$.  \hfill\qed

\medbreak

\proclaim {Lemma 5.3} Suppose that the $2$-dimensional subspace $X$ of
complex $\ell_\infty^k$ is spanned by two linearly independent vectors $y,z$
such that $|y|=|z|$. Then there is a linear isometry $\Phi$
of $l_\infty^k$ such that $\Phi y=\overline{\Phi z}$. In
particular, $\Phi(X)$ is spanned by two vectors $v,w$, all of
whose  coordinates are real and which satisfy
$(|v|^2+|w|^2)^{1/2}=|y|$.
\endproclaim

{\bf Proof.}  Write $y=(y_1,y_2,\cdots,y_k)$ and
$z=(z_1,z_2,\cdots,z_k)$. For
$j=1,2,\cdots,k$, let $\alpha_j$ be a complex number with
$|\alpha_j|=1$ such that $\alpha_j y_j= \overline{\alpha_j
z_j}$. Such numbers obviously exist, we may also impose
the condition $\Re\alpha_j y_j\ge0$. Now the isometry
$\Phi$ can be defined by the formula
$\Phi(x_1,x_2,\cdots,x_k)=(\alpha_1x_1,\alpha_2x_2,\cdots,\alpha_kx_k)$.
Clearly, the vectors $v={1\over2}(\Phi y+\Phi z)$ and
$w={1\over2i}(\Phi y-\Phi z)$ have the required property.
\hfill\qed
 \ms
Propositions 2.9, 2.10, and 5.1 suggest an alternate method for proving
Proposition 4.4 since they combine to take care of the case where the
operator achieves its norm at two ``flat" vectors:

\ms

\proclaim {Lemma 5.4} Let $X$ be a two-dimensional subspace of complex
$\ell_\infty^3$ and $T$ a complex-linear operator from $X$ into a Hilbert
space such that $\nm{Tx}=\nm{Ty}$, where $x,y$ are linearly
independent vectors in $X$ for which  $|x|=|y|$. Then $\pi_2(T)=\nm T$.
\endproclaim

{\bf Proof:} In view of Lemma 5.3 we can assume that there are vectors
$v$, $w$ in $X$ all of whose coordinates are real for which
$|v|^2+|w|^2=|x|^2$,   $x=v+iw$, and $y=v-iw$.  Thus if we let $E$ be the
collection of real-linear combinations of $\{v,w\}$, we can regard $X$
as the complexification $E_\C$ of $E$. The assumption on $\{x,y\}$ means
that the pair  $\{v,w\}$ satisfies  condition (iv) in Proposition 2.11, hence
condition (i) of Proposition 2.11 says that $T$ is the complexification of the
restriction of $T$ to $E$, whence by Proposition 2.10, $\pi_2(T)=\nm T$.
\hfill\qed

\medskip
The next lemma takes care of the case where the operator achieves its
norm at a non-flat vector.

\ms

\proclaim {Lemma 5.5} Suppose that $X$ is a $2$-dimensional subspace of
complex $\ell_\infty^3$ and  the norm one operator $X {\buildrel T \over
\longrightarrow} \ell_2^2$  achieves its norm at a non-flat vector
$x=(x_1,x_2,x_3)$ on the unit sphere of $X$.  Then there are norm one
operators $X {\buildrel V \over \longrightarrow} \ell_\infty^2$ and
$\ell_\infty^2 {\buildrel W \over \longrightarrow} \ell_2^2$ so that $T
= WV$.  Consequently, by Example 2.2,  $\pi_2(T) = 1$.
\endproclaim

{\bf Proof:} Since the result is trivial if $T$ has rank one, we assume
that $T$ has rank two.  This implies that two coordinates of $x$, say,
$x_1$ and $x_2$, are unimodular.  [Indeed, suppose, for example, that
$|x_2|$ and $|x_3|$ are both less than $1-\e$.  Take $y$ in $X$ with
$y_1=0$ and $0 < \nm{y} < \e$, so $\nm{x \pm y}=1$.  But since  $\ell_2$
is strictly convex, $\nm{T(x + \eta y)}> \nm{Tx}$ for either $\eta=1$ or
$\eta=-1$.]  We may also assume that $X$ contains two vectors, say $y,w$,
such that $y=(1,0,y_3)$ and $w=(0,1,w_3)$.  (Otherwise, $X$ is spanned by
two vectors with disjoint supports and the conclusion of the Lemma is
obvious.) Let
$\Psi$ be the function defined for $z\in\C$ by the formula
$$
\Psi(z)=\nm{T(x_1y+zw)}^2.
$$
Observe that $\Psi(z)=\nm{T(x_1y)+zTw}^2$ is a quadratic
function of $(\Re z,\Im z)$, which at infinity is
asymptotically equal to $m|z|^2$, where $m=\nm{Tw}^2>0$.
It follows that there is a number $z_0\in\C$ such that
$\Psi(z)=m|z-z_0|^2+\Psi(z_0)$ for every $z\in\C$. It is
obvious now that either $z_0=0$, so that $\Psi$ is
constant on the unit circle, or else $\Psi$ has a unique
local maximum on the unit circle (which must also be the
global maximum of $\Psi$ on the circle). Note that at
$z=x_2$ the function $\Psi$ does have a local maximum. In
each case it follows that, for every $z$ with $|z|\le1$,
we have $\Psi(z)\le\Psi(x_2)=1$.

Put $V(z_1y+z_2w)=(z_1,z_2)$. Then $V:X\rightarrow
l_\infty^2$ and $\nm V\le1$. The latter property of
$\Psi$ can be restated as follows: if $(z_1,z_2)\in\C^2$,
$|z_1|=|z_2|=1$, and $u=z_1/x_1$, then
$$
\nm{T(z_1y+z_2w)}=\sqrt{\Psi(z_2/u)}\le1=
\max\{|z_1|,|z_2|\}=\nm{V(z_1y+z_2w)}.
$$
Since $W=TV^{-1}$ attains its norm at an extreme point,
we have just checked that $\nm{W}\le1$.\hfill\qed
\ms
Now we can give an alternate:

\ms
{\bf Proof of Proposition 4.4:}  Suppose that the two-dimensional
subspace $X$ of complex $\ell_\infty^3$ fails the \twosump.  Let
$T_0$ be a norm one linear mapping of $X$ into a Hilbert space
$H$ whose $2$-summing norm is maximal among all norm one linear
maps of $X$ into $H$. Thus $\pi_2(T_0)>\nm{T_0}=1$; in
particular, $T_0$ is of rank $>1$.

By Lemma 5.5, $T_0$ does not attain its norm at any non-flat
vector. Hence, if $f_1,f_2\in X$ are the two flat vectors
described in Lemma 5.2, then
$\nm{T_0}=\max\{\nm{T_0f_1},\nm{T_0f_2}\}$.
Using Lemma 5.4, we rule out the possibility that
$\nm{T_0f_1}=\nm{T_0f_2}$.

Assume that $\nm{T_0}=\nm{T_0f_1}>\nm{T_0f_2}$.
Observe that for every $\epsilon>0$ there is an operator
$T_\epsilon:X\rightarrow H$ such that
$$
\nm{T_\epsilon-T_0}<\epsilon,\qquad
\nm{T_\epsilon f_1}=\nm{T_0f_1}
$$
and the inequality $\nm{T_\epsilon x}\le\nm{T_0x}$ is
possible only if $x=\lambda f_1$.

Since rank $T$ is $>1$, the latter property of
$T_\epsilon$ implies that $\pi_2(T_\epsilon)>\pi_2(T_0)$.
By the maximality of $\pi_2(T_0)$, we infer that
$\nm{T_\epsilon}>1$. However, $T_\epsilon$ does attain
its norm somewhere and  it cannot happen at any flat
vector, because as soon as
$\epsilon<\nm{T_0f_1}-\nm{T_0f_2}$ we have
$\nm{T_\epsilon f_2}<\nm{T_\epsilon f_1}=1$. Using
Lemma 5.5 again, we infer that
$\pi_2(T_\epsilon)=\nm{T_\epsilon}$. Now, letting
$\epsilon$ tend to $0$, we obtain that
$\pi_2(T_0)=\nm{T_0}$, which contradicts our initial
assumption. \hfill\qed

\ms
To find $3$-dimensional subspaces other than $\ell_\infty^3$
and $\ell_1^3$ which have the \twosump, it is natural to look
inside $\ell_\infty^4$.  However:

\ms

\proclaim {Proposition 5.6} Let $X$ be a three dimensional subspace of
complex $\ell_\infty^4$ not isometric to $\ell_\infty^3$. Then $X$
does not have the $2$-dimensional \twosump.
\endproclaim

{\bf Proof:} First notice that without loss of generality $X$ is
spanned by three vectors of the form
$(1,0,0,a_1),(0,1,0,a_2),(0,0,1,a_3)$ with $a_1,a_2,a_3$ non-negative
real numbers satisfying $a_1,a_2,a_3\le 1$ and $a_1+a_2+a_3>1$. Indeed,
let $(b_1,b_2,b_3,b_4)$ be a non-zero vector annihilating $X$ and
assume $|b_4|\ge |b_1|,|b_2|,|b_3|$. Applying $\ell_\infty^4$
isometries we may assume that $b_i/b_4$ are non-positive reals. Put
$a_i=-b_i/b_4$. Note that, since $X$ is not  isometric to
$\ell_\infty^3$, $a_1+a_2+a_3>1$.

Fix $\a,\b,\g,\d$ non-negative real numbers and $\varphi,\psi\in \C$
with $|\varphi|=|\psi|=1$ and consider the following two vectors in
$X$:
$$
\eqalign{ x=&(\a,\varphi,\g\psi,\a a_1+ \varphi a_2+\g \psi a_3)\cr
y=&(\b,0,-\d \psi,\b a_1-\d \psi a_3).}
$$
We are going to show that, for some choice of the parameters,
$|x|^2+|y|^2$ is constantly equal to one while for the same choice
$\span\{x,y\}$ does not contain a flat vector. Once this is proved one
concludes the proof as in Example 2.3.

$|x|^2+|y|^2\equiv 1$ is equivalent to
$$
1=\a^2+\b^2=\g^2+\d^2=|\a
a_1+\varphi a_2+\g \psi a_3|^2+|\b a_1-\d \psi a_3|^2 \leqno{(5.1)}
$$
while, if $|ax+by|\equiv 1$, we may assume without loss of generality
that $a=1$ and then
$$
1=|\a+b\b|=|\g -b\d |=|\a a_1+\varphi
a_2+\g \psi a_3+b( \b a_1-\d \psi a_3)|.\leqno{(5.2)}
$$
Note that if
$\b\not=0$ then the first equations in (5.1) and (5.2) imply that $2{\a
\over \b}\Re b=1-|b|^2$. Similarly, if $\d\not= 0$,
 $-2{\g \over \d}\Re b=1-|b|^2$. It follows that  $|b|=1$. If in
addition to $\b,\d>0$ also $\a>0$ or $\g>0$ then $b=\pm i$ and the last
equations in (5.1) and (5.2) imply that $\a a_1+\varphi a_2+\g \psi
a_3$ and
 $ \b a_1-\d \psi a_3$ are pointing in the same or opposite directions.
Thus, it is enough to find non-negative  reals $\a,\b,\g,\d$  with all
but possibly $\a$ or $\g$ positive and complex numbers
$\varphi$, $\psi$ of modulus one satisfying  (5.1) but
such that $\a a_1+\varphi a_2+\g \psi a_3$ and
 $ \b a_1-\d \psi a_3$ are not pointing in the same or opposite
directions.

Assume first that $a_2^2 + (a_1-a_3)^2<  1$ and $a_2>0$. Clearly there
are $0<\a,\g<1$  and $0<\varphi<\pi$ for which $|\a a_1+\varphi a_2+\g
a_3|>1$. Then also \hfill\break $|\a a_1+\varphi a_2+\g
a_3|^2+|(1-\a^2)^{1/2}a_1-(1-\g^2)^{1/2} a_3|^2>1$ . Replacing $\a,\g$
with $t\a,t\g$ for some $0<t<1$ we find $0<\a,\g<1$
 and $0<\varphi<\pi$ for which
$$
|\a a_1+\varphi a_2+\g  a_3|^2+|(1-\a^2)^{1/2}a_1-(1-\g^2)^{1/2}
a_3|^2=1.
$$
Clearly, $\a a_1+\varphi a_2+\g  a_3$ and
$(1-\a^2)^{1/2}a_1-(1-\g^2)^{1/2} a_3$ are not pointing in the same or
opposite direction.

If $a_{\pi(2)}^2 + (a_{\pi(1)}-a_{\pi(3)})^2\ge 1$ for all
permutations, $\pi$, of the indices $1,2,3$ for which $a_{\pi(2)}>0$
then, assuming as we may that $a_1\ge a_2,a_3$ and $a_1>0$, it is
easily checked that there are $\a,\b>0$ with $\a^2+\b^2=1$ for which
 $$
|\a a_1+a_2|^2+|\b a_1+a_3|^2>1>|\a a_1-a_2|^2+|\b a_1-a_3|^2.
$$
 Indeed,
$$
1\le (a_1^2 + (a_2-a_3)^2)^{1/2}<
a_1+(a_2^2+a_3^2)^{1/2}=\sup_{\a^2+\b^2=1}(|\a a_1+a_2|^2+|\b
a_1+a_3|^2)^{1/2}.
$$
Moreover, the $\sup$ is attained for $(\a,\b)$ proportional to
$(a_2,a_3)$.  For this choice of $(\a,\b)$,\break
 $ |\a a_1-a_2|^2+|\b a_1-a_3|^2 =( a_1 - (a_2^2+a_3^2)^{1/2})^2 < 1$.
 Choose now $\g=0, \d=1$ and notice that there is a one parameter
family of $\varphi,\psi$ for which $|\a a_1+\varphi a_2|^2+|\b a_1-\psi
a_3|^2=1$ but not for all of members of this family do $\a a_1+\varphi
a_2$ and $\b a_1-\psi a_3$ point in the same or opposite directions.
\hfill\qed

\ms

{\bf Remark.} It is easy to adjust the proof above to show that for
$n>3$ no $n$-dimensional subspace of $\ell_\infty^{n+1}$ has the
\twosump. Indeed without loss of generality any such subspace is
spanned by $n$ vectors of the form  $(1,0,\dots,0,a_1), \dots,$
$(0,\dots,0,1,a_n)$ with $0\le a_i\le 1$.
 If
$\sum a_i\le 1$ the subspace is isometric to $\ell_\infty^{n}$ which
does not have the \twosump. Otherwise these $n$ vectors can be blocked
to get three vectors which could replace the three vectors with which
we started the  proof above.

\ms

We next present a proof that all maximal distance subspaces of $L_1$ are
$\ell_1^n$ spaces. In particular there are no new subspaces of (real or
complex) $L_1$ with the \twosump . Essentially the same proof shows that
all maximal distance subspaces of $L_p$ are
$\ell_p^n$ spaces, $1\le p<\infty$. This fact is not new: it was observed by
J. Bourgain that the case $p=1$ follows from [FJ]. The case $1<p<2$ was
first proved in [BT]. Komorowski [Ko] was the first to prove the
$2<p<\infty$ case. The proof here is very similar to Komorowski's but includes
also the $1\le p<2$ case.

\proclaim{Proposition 5.7} Let $1\le
p<\infty$ and let $X$ be n-dimensional subspace of an $L_p(\mu)$ space with
$d(X,\ell_2^n) = n^{|{1\over p}-{1\over 2}|}$, then $X$ is isometric to
$\ell_p^n$.
\endproclaim

{\bf Proof.} By Lewis' theorem [L] we may assume that
$\mu$ is a probability measure and that $X$ has a basis
$x_1,x_2,\dots,x_n$ satisfying

$$\leqalignno{& \sum_{i=1}^n |x_i|^2 \equiv 1 &(5.3)\cr
\hbox{and}\quad \quad \ \ \ \ \ \ \ \ \ \ \ \  &\cr
& \int |\sum_{i=1}^n a_ix_i|^2  d\mu = {1\over n} \sum_{i=1}^n
|a_i|^2\ \ \quad\hbox{ for all scalars.}&(5.4)\cr}$$ Then, for $1\le p\le 2$,
$${1\over{n^{1/2}}}(\sum_{i=1}^n |a_i|^2)^{1/2}=(\int |\sum_{i=1}^n
a_ix_i|^2d\mu)^{1/2} \ge (\int |\sum_{i=1}^n a_ix_i|^pd\mu)^{1/p}$$
and
$$\leqalignno{{1\over{n^{1/2}}}(\sum_{i=1}^n |a_i|^2)^{1/2}=&(\int
|\sum_{i=1}^n a_ix_i|^2d\mu)^{1/2}\cr
\le &(\int |\sum_{i=1}^n a_ix_i|^pd\mu)^{1/2}\sup|\sum_{i=1}^n
a_ix_i|^{{2-p}\over 2}&(5.5)\cr
\le&(\int |\sum_{i=1}^n a_ix_i|^pd\mu)^{1/2}(\sum_{i=1}^n |a_i|^2)^{{2-p}\over
4} &(5.6)\cr}$$
(by (5.3)). Thus
$${1\over{n^{1/2}}}(\sum_{i=1}^n |a_i|^2)^{p/4}\le(\int
|\sum_{i=1}^n a_ix_i|^pd\mu)^{1/2}$$
and
$${1\over{n^{1/p}}}(\sum_{i=1}^n |a_i|^2)^{1/2}\le(\int
|\sum_{i=1}^n a_ix_i|^pd\mu)^{1/p}.$$
This shows that, if $T$ is the map sending the $x_i$'s to an orthonormal basis,
then \hfill\break $||T||\,||T^{-1}|| \le n^{|{1\over p}-{1\over 2}|}$.
It follows that there are $a_1,a_2,\dots,a_n$ for which equality is achieved
in both (5.5) and (5.6).

If $2<p<\infty$, then we get similarly that
$${1\over{n^{1/2}}}(\sum_{i=1}^n |a_i|^2)^{1/2} \le (\int |\sum_{i=1}^n
a_ix_i|^pd\mu)^{1/p}$$
and
$$\leqalignno{(\int |\sum_{i=1}^n a_ix_i|^pd\mu)^{1/p}
\le &(\int |\sum_{i=1}^n a_ix_i|^2d\mu)^{1/p}\sup|\sum_{i=1}^n
a_ix_i|^{{p-2}\over p}&(5.5')\cr
\le&(\int |\sum_{i=1}^n a_ix_i|^2d\mu)^{1/p}(\sum_{i=1}^n |a_i|^2)^{{p-2}\over
2p} &(5.6')\cr
=&{1\over{n^{1/p}}}(\sum_{i=1}^n |a_i|^2)^{1/2}.\cr}$$
Again we get that some $a_1,a_2,\dots,a_n$ must satisfy $(5.5')$ and
$(5.6')$ as equalities. Examining when equalities can occur in
$(5.5),(5.5'),(5.6)$  and
 $(5.6')$,
we see that for all $p$ there are $a_1,a_2,\dots,a_n$ such that
$|\sum_{i=1}^n a_ix_i|$ is a constant on its support, $A$, and the constant
must be $(\sum_{i=1}^n |a_i|^2)^{1/2}$ which we may assume is equal to $1$.
Moreover, on $A$,\hfill\break $(x_1(t),x_2(t),\dots,x_n(t))$ must be equal to
$\theta(t)(\bar a_1,\bar a_2, \dots,\bar a_n)$ for some function $\theta$
satisfying $|\theta(t)|\equiv 1$ . Applying a space isometry, we may assume
that $\theta\equiv 1$ and then  $\sum_{i=1}^n a_ix_i = \chi_{_A}$. Note also
that $\mu(A) = \int |\sum_{i=1}^n a_ix_i|^2d\mu = {1\over n}$. We thus get
that we may assume that $\chi_{_A} \in X$. Since each $x_i$ is constant on
$A$, we get that for all $f\in X$, $f_{|A}$ is a constant and $f_{|A^c}$
also belongs to $X$. Put
$$
Y=\{ f_{|A^c} \ :\ \ f\in X\};
$$
then $Y\subset X$ and dim$Y = n-1$. Necessarily $d(Y,\ell_p^{n-1}) =
(n-1)^{|{1\over p}-{1\over 2}|}$ and continue...\hfill$\qed$

\vfill\eject


\enddocument